\documentclass[12pt]{article}
\usepackage{amssymb}
\usepackage{amsthm}
\usepackage{amsmath}
\usepackage{graphicx}
 \newtheorem{thm}{Theorem}[section]

 \theoremstyle{definition}

 \numberwithin{equation}{section}

\begin{document}
\title{A Complete Global Solution to the Pressure Gradient Equation}
\author{Zhen Lei
\thanks{School of Mathematical Sciences, Fudan University,
Shanghai 200433, P. R. China; School of Mathematics and
Statistics, Northeast Normal University, Changchun 130024, P. R.
China; Key Laboratory of Mathematics for Nonlinear Sciences (Fudan
University), Ministry of Education, P. R. China, {\it email:
leizhn@yahoo.com.}} \and Yuxi Zheng\thanks{Department of
Mathematics, Pennsylvania State University, State college, PA
16802, {\it email: yzheng@math.psu.edu.}}}
\date{May 6, 2006}
\maketitle

\begin{abstract}
We study the domain of existence of a solution to a Riemann
problem for the pressure gradient equation in two space
dimensions. The Riemann problem is the expansion of a quadrant of
gas of constant state into the other three vacuum quadrants. The
global existence of a smooth solution was established in Dai and
Zhang [Arch.~Rational Mech.~Anal., {\bf 155}(2000), 277-298] up to
the free boundary of vacuum. We prove that the vacuum boundary
where the system is degenerate is the trivial coordinate axes.
\end{abstract}



\section{Introduction}

The pressure gradient system
\begin{equation}\label{1}
\left\{
\begin{array}{rcl}
u_t + p_x & = & 0, \\[2mm]
v_t + p_y & = & 0, \\[2mm]
E_t +(pu)_x + (pv)_y & = & 0,
\end{array}
\right.
\end{equation}
where $E = (u^2+ v^2)/2 + p$, appeared first in the flux-splitting
method of Li and Cao \cite{Licao} and Agarwal and Halt \cite{ah}
in numerical computation of the Euler system of a compressible
gas. Later, an asymptotic derivation was given in Zheng
\cite{zhengcpde, zhengyuxi}
 from the two-dimensional full Euler system for an ideal fluid
\begin{equation}\nonumber
\left\{
\begin{array}{rcl}
\rho_t + \nabla\cdot(\rho U) & = & 0, \\[2mm]
(\rho U)_t + \nabla\cdot(\rho U\otimes U + pI) & = & 0, \\[2mm]
(\rho E)_t + \nabla\cdot(\rho EU + pU) & = & 0,
\end{array}
\right.
\end{equation}
where $U = (u, v), E = (u^2+v^2)/2 + p/((\gamma-1)\rho)$, and
$\gamma>1$ is a gas constant. We refer the reader to the books of
Zheng \cite{zhengbook} and Li {\it et.~al.~}\cite{LiZhangYang} for
more background information, and the papers \cite{zhengcpde,
ZhangLiZhang, daizhang, zheng, Zheng, zhengyuxi} for recent
studies. After being decoupled from the pressure gradient system
\eqref{1}, the pressure satisfies the following second order
quasi-linear hyperbolic equation
\begin{equation}\label{1.1}
\Big(\frac{p_t}{p}\Big)_t - \Delta_{(x, y)}p = 0.
\end{equation}

Dai and Zhang \cite{daizhang} studied a Riemann problem for system
\eqref{1}, see also Yang and Zhang \cite{YangZhang} by the
hodograph method. In the self similar variables $\xi = x/t, \eta =
y/t$, the value of the pressure variable of the Riemann data is
\begin{equation}\label{1.3}
\begin{cases}
p(\xi, \eta) = \xi^2,  \quad \mbox{ for } 0 < \xi, \eta \leq
  \sqrt{p_1}, \quad (\xi - \sqrt{p_1})^2 + \eta^2 = p_1,\\[3mm]
p(\xi, \eta) = \eta^2, \quad \mbox{ for }  0 < \xi, \eta \leq
  \sqrt{p_1}, \quad \xi^2 + (\eta - \sqrt{p_1})^2 = p_1.
\end{cases}
\end{equation}
Here $p_1$ is any positive number. They showed that the Goursat
problem for system \eqref{1.1} admits a global solution in the
self-similar plane, which is smooth with a possible vacuum near
the origin (see Figure 1, where $a =\sqrt{p_1}$). {\sl We are
interested in the size of the vacuum boundary $\{(\xi, \eta)\ | \
p(\xi, \eta) = 0\}$ where the pressure gradient system \eqref{1.1}
is degenerate.} Somewhat surprisingly, our result shows that the
vacuum bubble is trivial and the entire vacuum boundary is the
trivial coordinate axes in the self-similar plane (see Figure 2,
where $a =\sqrt{p_1}$), which is stated in our main theorem at the
end of Section 3.

\begin{figure}
\medskip
\includegraphics [width=12cm,clip]{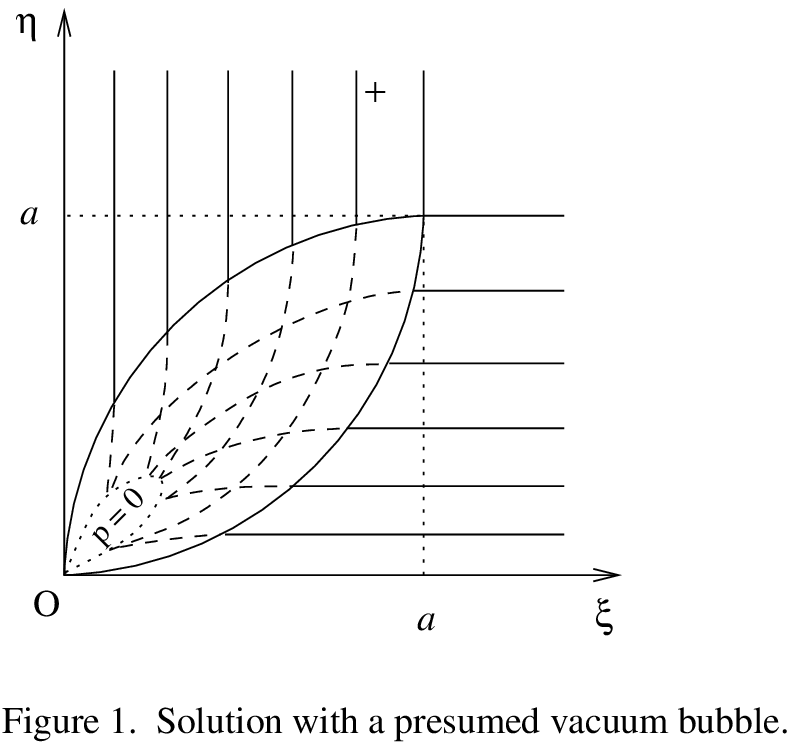}
\caption{}
\medskip
\end{figure}

Further motivation for the study of the current problem is that
the study of boundaries such as a sonic curve is important in
establishing the global existence of a solution to a general two
dimensional Riemann problem of the pressure gradient system. In
addition, the solution of the current problem covers wave
interaction problems in which only some fractions of the plane
waves are involved. Wave interactions of these kinds are common in
two dimensional Riemann problems. Finally, the study of the
pressure gradient system has motivated work on two dimensional
full Euler systems, see Li \cite{li} and Zheng
\cite{zheng04,zheng05}.

\begin{figure}
\medskip
\includegraphics [width=12cm,clip]{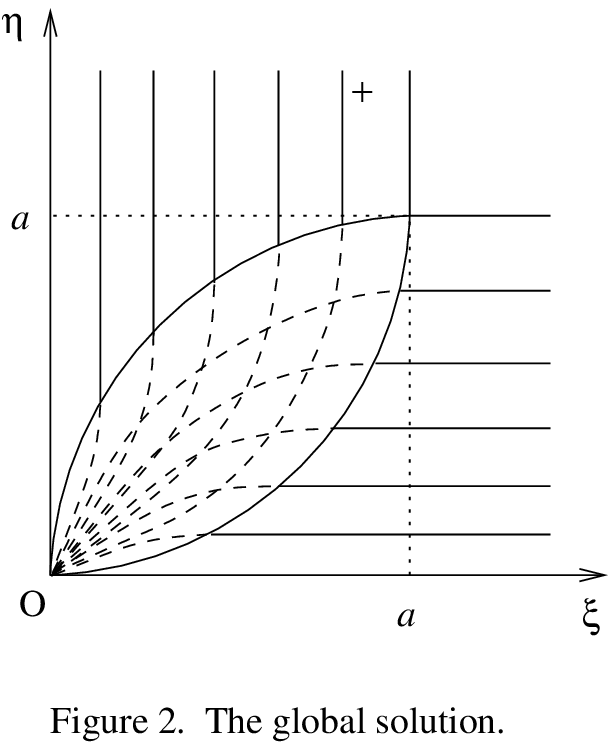}
\caption{}
\medskip
\end{figure}

\section{Integration along characteristics}

In the self-similar variables $\xi = x/t, \ \eta = y/t$, the
pressure gradient equation \eqref{1.1} takes the form
\begin{equation}\label{2.1}
\frac{(\xi\partial_\xi + \eta\partial_\eta)^2p}{p} - \Delta_{(\xi,
\eta)}p  + \frac{(\xi\partial_\xi + \eta\partial_\eta)p}{p} -
\frac{((\xi\partial_\xi + \eta\partial_\eta)p)^2}{p^2} = 0.
\end{equation}
In the polar coordinates
\begin{equation}\nonumber
\begin{cases}
r = \sqrt{\xi^2 + \eta^2},\\[3mm]
\theta = \arctan\frac{\eta}{\xi},
\end{cases}
\end{equation}
equation \eqref{2.1} can be decomposed along the characteristics
into the following form (see \cite{Zheng, LeiZheng})
\begin{equation}\label{z1}
\begin{cases}
\partial_+\partial_-p =  mp_r\,\partial_-p,\\[3mm]
\partial_-\partial_+p =  - mp_r\,\partial_+p,
\end{cases}
\end{equation}
provided that
$$p < r^2,$$
where
$$m = \frac{\lambda r^4}{2p^2},$$
and
\begin{equation}\label{2.3}
\begin{cases}
\partial_\pm = \partial_\theta \pm \frac{1}{\lambda}\partial_r,\\[3mm]
\lambda = \sqrt{\frac{p}{r^2(r^2 - p)}}.
\end{cases}
\end{equation}

Note that the equation is invariant under the following scaling
transformation
$$
(\xi, \eta, p) \longrightarrow (\frac{\xi}{\sqrt{p_1}},
\frac{\eta}{\sqrt{p_1}}, \frac{p}{p_1}). \qquad (p_1 > 0)
$$
Thus, without loss of generality, the corresponding boundary
condition \eqref{1.3} in the polar coordinates can be set as the
form
\begin{equation}\label{a1}
\begin{cases}
p = \xi^2 = r^2\cos^2\theta = 4\cos^4\theta, \quad {\rm on}\quad
  r = 2\cos\theta, \quad \pi/4 \leq \theta \leq \pi/2; \\[3mm]
p = \eta^2 = r^2\sin^2\theta = 4\sin^4\theta, \quad{\rm on}\quad r
  =2\sin\theta, \quad 0 \leq \theta \leq \pi/4.
\end{cases}
\end{equation}
The solution exists in the interaction zone up to a possible
vacuum bubble. See Figure 3.

\begin{figure}
\medskip
\includegraphics [width=12cm,clip]{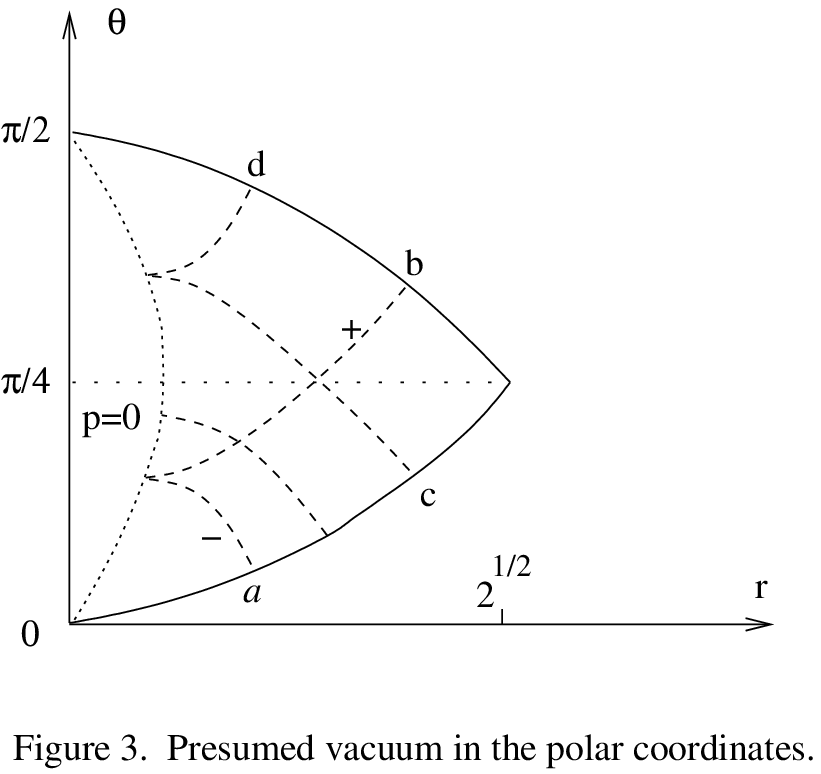}
\caption{}
\medskip
\end{figure}

The characteristic form \eqref{z1} of the pressure gradient
equation enjoys a number of useful properties. For example, the
quantities $\partial_\pm p$ keep their positivities/negativities
along characteristics of a plus/minus family, and the
sign-persevering quantities yield monotonicity of the primary
variable $p$ (see \cite{Zheng}, also \cite{LeiZheng} where the
authors propose to call them {\it Riemann sign-persevering
variables}), and the fact that a state adjacent to a constant
state for the pressure gradient system must be a simple wave in
which $p$ is constant along the characteristics of a plus/minus
family \cite{Zheng}.

The characteristic decomposition \eqref{z1} has played an
important role and was a powerful tool for building the existence
of smooth solutions in the work of Dai and Zhang \cite{daizhang}.
We are interested in the size of the vacuum boundary $\{(r,
\theta) \ |\ p(r, \theta) = 0\}$ where the pressure gradient
equation is degenerate. For this purpose, we rewrite \eqref{z1} as
\begin{equation}\label{2.2}
\begin{cases}
\partial_+\partial_-p =  q\partial_+p\partial_-p
  - q(\partial_-p)^2,\\[3mm]
\partial_-\partial_+p = q\partial_+p\partial_-p
  - q(\partial_+p)^2,
\end{cases}
\end{equation}
where
\begin{equation}\label{2.4}
q = \frac{r^2}{4p(r^2 - p)}.
\end{equation}
Define the characteristic curves $r_-^a(\theta)$ and
$r_+^b(\theta)$ as
\begin{equation}\label{2.5}
\begin{cases}
\frac{dr_-^a(\theta)}{d\theta} = - \frac{1}{\lambda
  \big(r_-^a(\theta), \theta\big)},\\[3mm]
r_-^a(\theta_a) = 2\sin\theta_a,
\end{cases}
{\rm and}\quad
\begin{cases}
\frac{dr_+^b(\theta)}{d\theta} = \frac{1}{\lambda
  \big(r_+^b(\theta), \theta\big)},\\[3mm]
r_+^b(\theta_b) = 2\cos\theta_b.
\end{cases}
\end{equation}
We point out here that for convenience, we also use the notation
$r_-^a(r, \theta)$ ($r_+^b(r, \theta)$, respectively) which
represents the characteristic passing through the point $(r,
\theta)$ and intersecting the lower (upper, respectively) boundary
at point $a$ ($b$, respectively). See Figure 3.

Now let us rewrite the system \eqref{2.2} as the following form
\begin{equation}\nonumber
\begin{cases}
\partial_+\Big(\frac{1}{\partial_-p}\exp\int_{\theta_b}^\theta
  q\partial_+p\big(r_+^b(\phi), \phi\big)d\phi\Big)
  = q\exp\int_{\theta_b}^\theta q\partial_+p
  \big(r_+^b(\phi), \phi\big)d\phi,\\[3mm]
\partial_-\Big(\frac{1}{\partial_+p}\exp\int_{\theta_a}^\theta
  q\partial_-p\big(r_-^a(\phi), \phi\big)d\phi\Big)
  = q\exp\int_{\theta_a}^\theta q\partial_-p
  \big(r_-^a(\phi), \phi\big)d\phi.
\end{cases}
\end{equation}
Integrate the above equations along the positive and negative
characteristics $r_+^b(\theta)$ and $r_-^a(\theta)$ from
$\theta_b$ and $\theta_a$ to $\theta$, respectively, with respect
to $\theta$, one can get the iterative expressions of
$\partial_+p$ and $\partial_-p$:
\begin{equation}\label{2.6}
\begin{cases}
\frac{1}{\partial_-p}\exp\int_{\theta_b}^\theta
  q\partial_+p\big(r_+^b(\phi), \phi\big)d\phi = \frac{1}{\partial_-p}
  \big(2\cos\theta_b, \theta_b\big)\\[3mm]
\quad +\ \int_{\theta_b}^\theta q\big(r_+^b(\psi),
  \psi\big)\exp\int_{\theta_b}^{\psi} q
  \partial_+p\big(r_+^b(\phi), \phi\big)d\phi d\psi,\\[3mm]
\frac{1}{\partial_+p}\exp\int_{\theta_a}^\theta
q\partial_-p\big(r_-^a(\phi),
  \phi\big)d\phi = \frac{1}{\partial_+p}\big(2\sin\theta_a,
\theta_a\big)\\[3mm]
\quad +\ \int_{\theta_a}^\theta q\big(r_-^a
  (\psi), \psi\big)\exp\int_{\theta_a}^{\psi}
  q\partial_-p\big(r_-^a(\phi), \phi\big)d\phi d\psi.
\end{cases}
\end{equation}

On the other hand, noting the boundary condition \eqref{a1}, a
straight forward calculation shows that
\begin{eqnarray}\label{2.7}
&&\exp\int_{\theta_b}^\theta q\partial_+p\big(r_+^b(\phi),
  \phi\big)d\phi\\\nonumber
&&= \exp\frac{1}{4}\int_{\theta_b}^\theta \Big(\frac{1}{p} +
  \frac{1}{r^2 - p}\Big)\partial_+p\big(r_+^b(\phi),
  \phi\big)d\phi\\\nonumber
&&= \frac{p^{\frac{1}{4}}\big(r_+^b(\theta), \theta\big)}
  {\sqrt{2}\cos\theta_b}\exp\Big\{\frac{1}{4}\int_{p(2\cos\theta_b,
  \theta_b)}^{p\big(r_+^b(\theta), \theta\big)}
  \frac{1}{r^2(r_+^b) - p_1}dp_1\Big\}.
\end{eqnarray}
Similarly, one has
\begin{eqnarray}\label{2.8}
&&\exp\int_{\theta_a}^\theta q\partial_-p\big(r_-^a(\phi),
  \phi\big)d\phi\\\nonumber
&&= \frac{p^{\frac{1}{4}}\big(r_-^a(\theta), \theta\big)}
  {\sqrt{2}\sin\theta_a}\exp\Big\{\frac{1}{4}\int_{p(2\sin\theta_a,
  \theta_a)}^{p\big(r_-^a(\theta), \theta\big)}
  \frac{1}{r^2(r_-^a) - p_1}dp_1\Big\}.
\end{eqnarray}
Then, we have
\begin{eqnarray}\label{2.9}
&&\int_{\theta_b}^\theta q\big(r_+^b(\psi),
  \psi\big)\Big\{\exp\int_{\theta_b}^{\psi} q
  \partial_+p\big(r_+^b(\phi), \phi\big)d\phi\Big\}
  d\psi\\\nonumber
&&= \frac{1}{4\sqrt{2}\cos\theta_b}\int_{\theta_b}^
  \theta\frac{r^2p^{- \frac{3}{4}}}{r^2 - p}\big(r_+^b(\psi),
  \psi\big)\\\nonumber
&&\quad \times\exp\Big\{\frac{1}{4}\int_{p(2\cos\theta_b,
  \theta_b)}^{p\big(r_+^b(\psi), \psi\big)}
  \frac{1}{r^2(r_+^b) - p_1}dp_1\Big\}d\psi,
\end{eqnarray}
and
\begin{eqnarray}\label{2.10}
&&\int_{\theta_a}^\theta q\big(r_-^a
  (\psi), \psi\big)\exp\int_{\theta_a}^{\psi}
  q\partial_-p\big(r_-^a(\phi), \phi\big)d\phi d\psi\\\nonumber
&&= \frac{1}{4\sqrt{2}\sin\theta_a}\int_{\theta_a}^
  \theta\frac{r^2p^{- \frac{3}{4}}}{r^2 - p}\big(r_-^a(\psi),
\psi\big)\\\nonumber &&\quad
\times\exp\Big\{\frac{1}{4}\int_{p(2\sin\theta_a,
  \theta_a)}^{p\big(r_-^a(\psi), \psi\big)}\frac{1}{r^2(r_-^a) -
p_1}dp_1\Big\} d\psi.
\end{eqnarray}

Finally, by substituting \eqref{2.8} and \eqref{2.10} into the
second equality of \eqref{2.6}, we arrive at a new iterative
expression for $\partial_+p$ :
\begin{eqnarray}\label{2.11}
&&\partial_+p\big(r, \theta\big)\\\nonumber &=&
  \frac{\exp\int_{\theta_a}^\theta q\partial_-p\big(r_-^a(\phi),
  \phi\big)d\phi}{\frac{1}{\partial_+p}\big(2\sin\theta_a,
  \theta_a\big) + \int_{\theta_a}^\theta q\big(r_-^a
  (\psi), \psi\big)\exp\int_{\theta_a}^{\psi}
  q\partial_-p\big(r_-^a(\phi), \phi\big)d\phi d\psi}\\\nonumber
&=& \frac{4p^{\frac{1}{4}}(r, \theta)
  \exp\Big\{\frac{1}{4}\int_{p(2\sin\theta_a, \theta_a)}^
  {p\big(r_-^a(\theta), \theta\big)}\frac{1}{r^2(r_-^a) - p_1}
  dp_1\Big\}}{\frac{1}{2\sqrt{2}\sin^2\theta_a\cos\theta_a}
  + \int_{\theta_a}^\theta\frac{r^2p^{- \frac{3}{4}}}{r^2 - p}
  \big(r_-^a(\psi), \psi\big)\exp\Big\{\frac{1}{4}
  \int_{p(2\sin\theta_a, \theta_a)}^{p\big(r_-^a(\psi),
  \psi\big)}\frac{1}{r^2(r_-^a) - p_1}dp_1\Big\} d\psi}.
\end{eqnarray}
Similarly, by substituting \eqref{2.7} and \eqref{2.9} into the
first equality of \eqref{2.6}, one obtains the new iterative
formula for $\partial_-p$ :
\begin{eqnarray}\label{2.12}
&&- \partial_-p(r, \theta)\\\nonumber &=&
  \frac{- \exp\int_{\theta_b}^\theta q\partial_+p\big(r_+^b(\phi),
  \phi\big)d\phi}{\frac{1}{\partial_-p}\big(2\cos\theta_b,
  \theta_b\big) + \int_{\theta_b}^\theta q\big(r_+^b
  (\psi), \psi\big)\exp\int_{\theta_b}^{\psi}
  q\partial_+p\big(r_+^b(\phi), \phi\big)d\phi d\psi}\\\nonumber
&=& \frac{4p^{\frac{1}{4}}(r, \theta)
  \exp\Big\{\frac{1}{4}\int_{p(2\cos\theta_b,
  \theta_b)}^{p\big(r_+^b(\theta), \theta\big)}
  \frac{1}{r^2(r_+^b) - p_1}dp_1\Big\}}{\frac{1}{2\sqrt{2}\cos^2
  \theta_b\sin\theta_b} - \int_{\theta_b}^
  {\theta}\frac{r^2p^{- \frac{3}{4}}}{r^2 - p}\big(r_+^b
  (\psi), \psi\big)\exp\Big\{\frac{1}{4}\int_{p(2\cos\theta_b,
  \theta_b)}^{p\big(r_+^b(\psi), \psi\big)}\frac{1}{r^2(r_+^b)
  - p_1}dp_1\Big\} d\psi}.
\end{eqnarray}

\section{The boundary of the  vacuum bubble}

In this section, we use the iterative formulas \eqref{2.11} and
\eqref{2.12} to prove that the vacuum bubble $\{(r, \theta)\ |\
p(r, \theta) = 0, \theta\in(0, \pi/2)\}$ is in fact the trivial
origin $\{(0, 0)\}$ in the self-similar plane. We use the method
of contradiction. Assume to the contrary that there is a bubble
with boundary $r_0(\theta)\ge 0$ for all $\theta \in (0,
\frac{\pi}{2})$, and $r_0(\theta)> 0$ for some $\theta \in (0,
\frac{\pi}{2})$. That means $p\big(r_0(\theta), \theta\big) = 0$
for $\theta \in (0, \frac{\pi}{2})$, and the solution is smooth in
the domain bounded by the bubble $r_0(\theta)$ and the upper and
lower characteristic boundaries. We intend to deduce
contradictions, which in turn proves that $r_0(\theta) = 0$ for
all $\theta \in (0, \pi/2)$ and the vacuum bubble $\{(r_0(\theta),
\theta)\}$ is in fact the trivial origin $\{(0, 0)\}$ in the
self-similar plane.

Before we start the above procedure, we point out two observations
which roughly imply the nonexistence of the vacuum bubble. For
presenting the observations, we assume further that $r_0(\theta) >
0$ for all $\theta \in (0, \frac{\pi}{2})$.

First, we can compute easily the characteristic slope
$$
d\theta/dr = -\lambda = -\frac{1}{2\sin\theta} \sim - \frac{1}{2}
$$
along the upper characteristic boundary and consider $\theta\to
\pi/2$. Similar result holds on the lower boundary. Now assuming
there is a bubble, noting that along the plus characteristics, by
definition of $\lambda$ in \eqref{2.3}, we see easily that
$$
d\theta/dr = \lambda = \sqrt{\frac{p}{r^2(r^2 - p)}} \rightarrow 0
$$
as $p \rightarrow 0$ (near the vacuum bubble) except for $(0, 0)$
and $(0, \frac{\pi}{2})$. The above computation reveals that there
is some kind of inconsistency for the slopes of the
characteristics at $(0, 0)$ and $(0, \frac{\pi}{2})$ in the polar
coordinate plane.

Next, let us calculate the decay rate of $p$ along the middle line
$\theta=\pi/4$. By using the symmetry of system \eqref{2.1} (see
Figure 3), we can obtain that
$$
\partial_\pm p \sim\pm M_0p^{1/2}
$$
asymptotically. To show the details, we propose
$$
\partial_\pm p \sim \pm M_0p^{\frac{1}{4} + \delta},
$$
as $(r, \theta)$ tends to a point of the vacuum bubble. Then, by
\eqref{2.11}, we have
\begin{eqnarray}\nonumber
\partial_+p\big(r, \theta\big) &=& 4p^{\frac{1}{4}}(r,
  \theta)\exp\Big\{\frac{1}{4}\int_{p(2\sin\theta_a, \theta_a)}^
  {p\big(r_-^a(\theta), \theta\big)}\frac{1}{r^2(r_-^a) - p_1}
dp_1\Big\}\Big/\Big\{\frac{1}{2\sqrt{2}\sin^2\theta_a\cos\theta_a}\\\nonumber
& & + \int_{\theta_a}^\theta\frac{r^2p^{- \frac{3}{4}}}{(r^2 -
  p)\partial_-p}\big(r_-^a(\psi), \psi\big)\exp\Big\{\frac{1}{4}
  \int_{p(2\sin\theta_a, \theta_a)}^{p\big(r_-^a(\psi),
  \psi\big)}\frac{1}{r^2(r_-^a) - p_1}dp_1\Big\} dp\Big\}\\\nonumber
&=& \frac{M_0}{4\delta}p^{\frac{1}{4} + \delta} + {\rm high\
order\ terms},
\end{eqnarray}
which implies that $\delta$ should be $\frac{1}{4}$.

Thus, by using \eqref{2.3}, we have
$$
\partial_r p \sim \frac{M_0}{r^2}p
$$
asymptotically as $(r, \theta)$ tends to a point of the vacuum
bubble. Therefore,
\begin{equation}\label{expdec}
p \sim c \,\exp\left({-\frac{M_0}{r}}\right)
\end{equation}
asymptotically as $(r, \theta)$ tends to a point of the vacuum
bubble on the line $\theta = \pi/4$, which implies that there is
no interior vacuum at least at $\theta = \frac{\pi}{4}$.

We point out incidentally that Zheng's previous numerical
computation of the bubble, referred to in Dai and Zhang
\cite{daizhang},
 is probably caused by the fast exponential decay
\eqref{expdec}.

\begin{figure}
\medskip
\includegraphics [width=12cm,clip]{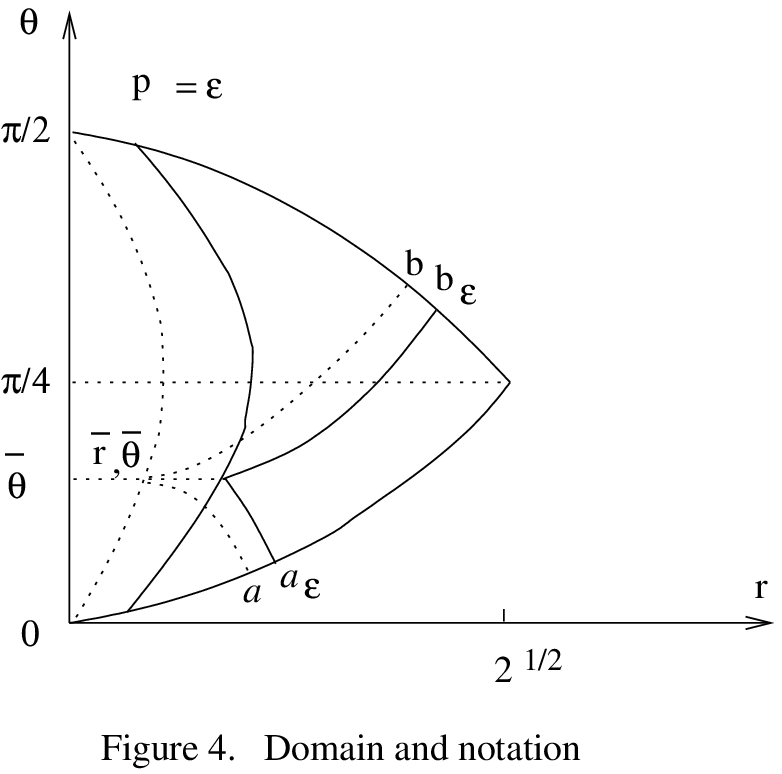}
\caption{}
\medskip
\end{figure}

In what follows, we concentrate on establishing the above formal
intuition rigorously. With the symmetry of system \eqref{2.2}, let
us restrict our arguments on $\theta \in (0, \frac{\pi}{4}]$.

Fix a  point $(\overline{r}, \overline{\theta})$ on the bubble
boundary. Let us denote $\mathbb{D}(\overline{r},
\overline{\theta})$ the bounded domain determined by the positive
and negative characteristic curves starting from $(\overline{r},
\overline{\theta})$ and $(\sqrt{2}, \frac{\pi}{4})$. See Figure 4.

For $0 < \epsilon < 1$, define a curve $r_\epsilon(\theta)$,
$\theta \in (0, \frac{\pi}{2})$, by
\begin{equation}\label{3.1}
p\big(r_\epsilon(\theta), \theta\big) = \epsilon.
\end{equation}
From \eqref{2.11} and \eqref{2.12}, it is easy to see that
\begin{equation}\label{3.2}
\begin{cases}
\partial_+p(r, \theta) > 0,\\
\partial_-p(r, \theta) < 0,
\end{cases}
\quad {\rm for}\quad 0 < \theta \leq \frac{\pi}{4}.
\end{equation}
Thus, by \eqref{2.3} and \eqref{3.2}, we have
\begin{equation}\label{3.3}
p_r(r, \theta) > 0,
\end{equation}
which implies that the curve $r_\epsilon(\theta)$,  $\theta \in
(0, \frac{\pi}{2})$, defined in \eqref{3.1} is smooth if $\epsilon
\in (0, 1)$. Here we point out that sometimes we still denote the
characteristics passing through $\big(r_\epsilon(\theta),
\theta\big)$ or any other point $(r, \theta)$ by $r_-^a(\theta)$
and $r_+^b(\theta)$ or $r_-^a(r, \theta)$ and $r_+^b(r, \theta)$,
when in fact $a = a_\epsilon$ is dependent on $\epsilon$, etc. We
caution the reader that the intersection is determined by the
point $(r, \theta)$ and the characteristic, which is hidden for
notational convenience and always well understood.

Now we fix $\epsilon_0 = 1$ ($\epsilon_0$ can be any fixed
positive constant in $(0, 1]$, which do not change any of the
following arguments). Define
\begin{equation}
M_1 = \max_{\mbox{ over } S} \{p^{- \frac{1}{2}}\partial_+p, \ -
p^{- \frac{1}{2}}\partial_-p\},
\end{equation}
where
$$
S :=\{ (r, \theta)\ | \ \theta\in (0, \pi/2), \
r_{\epsilon_0}(\theta) \leq r \leq \sqrt{2}\}.$$ Then we have
\begin{equation}\label{z2}
\begin{cases}
\begin{array}{rcl}
\partial_+p  & \leq & M_1p^{\frac{1}{2}},\\
- \partial_-p & \leq & M_1p^{\frac{1}{2}},
\end{array}
\end{cases}
\end{equation}
for all $(r, \theta)$ with $r \geq r_{\epsilon_0}(\theta)$. Note
that $M_1$ depends only on $\epsilon_0$ and does not depend on
$(\overline{r}, \overline{\theta})$.

Next, let
\begin{equation}\label{z4}
\begin{cases}
A(r, \theta) = \frac{\exp\Big\{- \frac{1}{4}\int_{p(2\sin\theta_a,
\theta_a)}^
  {p\big(r_-^a(\theta), \theta\big)}\frac{1}{r^2(r_-^a) - p_1}
  dp_1\Big\}}{2\sqrt{2}\sin^2\theta_a\cos\theta_a},\\[4mm]
B(r, \theta) =
  \frac{\exp\Big\{-\frac{1}{4}\int_{p(2\cos\theta_b,
  \theta_b)}^{p\big(r_+^b(\theta), \theta\big)}
  \frac{1}{r^2(r_+^b) - p_1}dp_1\Big\}}{2\sqrt{2}\cos^2
  \theta_b\sin\theta_b}.
\end{cases}
\end{equation}
Then, let
\begin{eqnarray}\label{z6a}
M_2 = \max\big\{\max_{(r, \theta) \in \mathbb{D}(\overline{r},
\overline{\theta})}\frac{4p^{- \frac{1}{4}}(2\sin\theta_a,
\theta_a)}{A(r, \theta)}, \max_{(r, \theta) \in
\mathbb{D}(\overline{r}, \overline{\theta})}\frac{4p^{-
\frac{1}{4}}(2\cos\theta_b, \theta_b)}{B(r, \theta)}\big\},
\end{eqnarray}
and
\begin{eqnarray}\label{z6}
M_3 = \max\big\{M_1, M_2 + 1\big\}.
\end{eqnarray}
Just as what we had pointed out before, the intersections $a$ and
$b$ in \eqref{z4} and \eqref{z6a} vary with $(r, \theta)$ and
characteristics, which are not expressed explicitly for notational
convenience. We intend to prove that the inequalities \eqref{z2}
are still valid for all points $(r, \theta) \in
\mathbb{D}(\overline{r}, \overline{\theta})$ with $M_1$ being
replaced by $M_3$. Namely, there hold
\begin{equation}\label{z3}
\begin{cases}
\begin{array}{rcl}
\partial_+p   & \leq & M_3p^{\frac{1}{2}},\\
- \partial_-p & \leq & M_3p^{\frac{1}{2}},
\end{array}
\end{cases}
(r, \theta) \in \mathbb{D}(\overline{r}, \overline{\theta}).
\end{equation}
Note that the positive constant $M_3 (\geq M_1)$ is independent of
$(r, \theta) \in \mathbb{D}(\overline{r}, \overline{\theta})$, but
depends on the fixed $(\overline{r}, \overline{\theta})$.

We start to prove  \eqref{z3}. Suppose that \eqref{z3} is correct
up to a line segment $r_\epsilon(\theta) \subset \mathbb{D}
(\overline{r}, \overline{\theta})$, then we improve \eqref{z3} to
strict inequalities on the line segment $r_\epsilon(\theta)$ in
the same domain. In fact, for $(r, \theta) \in
\mathbb{D}(\overline{r}, \overline{\theta})$ with $r \leq
r_{\epsilon_0}(\theta)$, let us compute
\begin{eqnarray}\label{z7}
\partial_+p(r, \theta) &=&  4p^{\frac{1}{4}}(r, \theta)
  \Big/\Big\{A(r, \theta)
  + \int_{\theta_a}^\theta\frac{r^2p^{- \frac{3}{4}}}{r^2 -
  p}\big(r_-^a(\psi), \psi\big)\\\nonumber
&&\times \exp\Big\{\frac{1}{4}
  \int_{p\big(r_-^a(\theta), \theta\big)}^{p\big(r_-^a(\psi),
  \psi\big)}\frac{1}{r^2(r_-^a) - p_1}dp_1\Big\}
  d\psi\Big\}\\\nonumber
&=& 4p^{\frac{1}{4}}(r, \theta)\Big/\Big\{
  A(r, \theta) + \int_{\theta_a}^\theta\frac{-
  \partial_-pr^2p^{- \frac{3}{4}}}{(-\partial_-p)(r^2 - p)}
  \big(r_-^a(\psi), \psi\big)\\\nonumber
& & \times \exp\Big\{\frac{1}{4}
  \int_{p\big(r_-^a(\theta), \theta\big)}^{p\big(r_-^a(\psi),
  \psi\big)}\frac{1}{r^2(r_-^a) - p_1}dp_1\Big\}d\psi\Big\}\\\nonumber
&\leq& M_3p^{\frac{1}{4}}(r, \theta)\Big/\Big\{
  \frac{M_3A(r, \theta)}{4} +
  \int_{p(2\sin\theta_a, \theta_a)}^
  {p\big(r_-^a(\theta), \theta\big)}
  \frac{- \frac{1}{4}p^{- \frac{5}{4}}r^2}{r^2 - p}
  \big(r_-^a(\psi), \psi\big)\\\nonumber
& & \times \exp\Big\{\frac{1}{4}
  \int_{p\big(r_-^a(\theta), \theta\big)}^{p\big(r_-^a(\psi),
  \psi\big)}\frac{1}{r^2(r_-^a) - p_1}dp_1\Big\}dp\Big\}\\\nonumber
&=& \frac{M_3p^{\frac{1}{4}}(r, \theta)
  }{\frac{M_3A(r, \theta)}{4} +
f(\theta^-)\int_{p(2\sin\theta_a, \theta_a)}^
  {p\big(r_-^a(\theta), \theta\big)} - \frac{1}{4}p^{-
  \frac{5}{4}}dp}\\\nonumber
&=& \frac{M_3p^{\frac{1}{4}}(r, \theta)
  }{\frac{M_3A(r, \theta)}{4} +
f(\theta^-)\big[p^{- \frac{1}{4}}(r, \theta) -
  p^{- \frac{1}{4}}(2\sin\theta_a, \theta_a)\big]}\ ,
\end{eqnarray}
where
\begin{eqnarray}\label{z8}
f(\theta^-) &=& \frac{r^2} {r^2 - p}\big(r_-^a(\theta^-),
  \theta^-\big)\exp\Big\{\frac{1}{4}\int_{p\big(r_-^a(\theta),
  \theta\big)}^{p\big(r_-^a(\theta^-), \theta^-\big)}
  \frac{1}{r^2(r_-^a) - p_1}dp_1\Big\}\\\nonumber
&>& \exp\Big\{\frac{1}{4}\int_{p\big(r_-^a(\theta),
  \theta\big)}^{p\big(r_-^a(\theta^-), \theta^-\big)}
  \frac{1}{\sqrt{2}^2}dp_1\Big\}\\\nonumber
&=& \exp\Big\{\frac{p\big(r_-^a(\theta^-),
  \theta^-\big) - p\big(r_-^a(\theta),
  \theta\big)}{8}\Big\} > 1
\end{eqnarray}
with some $\theta_a < \theta^- < \theta$ and $A(r, \theta)$ is
defined in \eqref{z4}. Thus, by \eqref{z6}, \eqref{z7} and
\eqref{z8}, we have
\begin{eqnarray}\label{3.9}
\partial_+p(r, \theta) &<&
\frac{M_3p^{\frac{1}{4}}(r, \theta)
  }{\frac{M_3A(r, \theta)}{4} + \big[p^{- \frac{1}{4}}(r, \theta) -
  p^{- \frac{1}{4}}(2\sin\theta_a, \theta_a)\big]}\\\nonumber
&\leq& M_3p^{\frac{1}{2}}(r, \theta).
\end{eqnarray}

Similarly, by \eqref{2.6}, \eqref{2.7} and \eqref{2.9}, we have
\begin{eqnarray}
- \partial_-p(r, \theta) &=&
  4p^{\frac{1}{4}}(r, \theta)\Big/\Big\{B(r,
  \theta) - \int_{\theta_b}^{\theta}\frac{r^2p^{- \frac{3}{4}}}{r^2
  - p}\big(r_+^b(\psi), \psi\big)\\\nonumber
&\times& \exp\Big\{\frac{1}{4}\int_{p\big
  (r_+^b(\theta), \theta\big)}^{p\big(r_+^b(\psi), \psi\big)}
  \frac{1}{r^2(r_+^b) - p_1}dp_1\Big\} d\psi\Big\}\\\nonumber
&\leq& \frac{M_3p^{\frac{1}{4}}(r, \theta)
  }{\frac{M_3B(r, \theta)}{4} + g(\theta^+)\big[p^{- \frac{1}{4}}(r,
\theta) -
  p^{- \frac{1}{4}}(2\cos\theta_b, \theta_b)\big]}\ ,
\end{eqnarray}
where
\begin{eqnarray}\nonumber
g(\theta^+) &=& \frac{r^2}{r^2 - p}\big(r_+^b(\theta^+),
  \theta^+\big)\exp\Big\{\frac{1}{4}\int_{p\big(r_+^b(\theta),
  \theta\big)}^{p\big(r_+^b(\theta^+), \theta^+\big)}
  \frac{1}{r^2(r_+^b) - p_1}dp_1\Big\}\\\nonumber
&>& \exp\Big\{\frac{1}{4}\int_{p\big
  (r_+^b(\theta), \theta\big)}^{p\big(r_+^b(\theta^+), \theta^+\big)}
  \frac{1}{\sqrt{2}^2}dp_1\Big\}\\\nonumber
&=& \exp\Big\{\frac{p\big(r_+^b(\theta^+), \theta^+\big)
  - p\big(r_+^b(\theta), \theta\big)}{8}\Big\} > 1
\end{eqnarray}
with some $\theta < \theta^+ < \theta_b$ and $B(r, \theta)$ is
given in \eqref{z4}. Thus, the similar estimate as \eqref{3.9}
holds:
\begin{eqnarray}\label{3.13}
- \partial_-p(r, \theta) < M_3p^{\frac{1}{2}}(r, \theta).
\end{eqnarray}
Since $M_3$ depends only on $(\overline{r}, \overline{\theta})$,
particularly, is independent of $(r, \theta) \in
\mathbb{D}(\overline{r}, \overline{\theta})$, the proof of
\eqref{z3} follows by \eqref{3.9} and \eqref{3.13}.

Now with the aid of \eqref{2.3}, we add up \eqref{3.9} and
\eqref{3.13} to yield
\begin{eqnarray}\nonumber
p_r \leq \frac{M_3}{\sqrt{r^2(r^2 - p)}}p \leq \frac{2
M_3}{\overline{r}^2}p,
\end{eqnarray}
for all $(r, \theta)$ in a small neighborhood of
$(\overline{r},\overline{\theta})$ in $\mathbb{D}(\overline{r},
\overline{\theta})$. Thus, a simple integration of the above
inequality with respect to $r$ from $\overline{r}$ to $r$ yields
\begin{eqnarray}\label{z5}
p(\overline{r}, \overline{\theta}) \geq p(r,
\overline{\theta})\exp\{- \frac{M_3}{\overline{r}^2}(r -
\overline{r})\} \quad {\rm for}\quad \overline{r} < r.
\end{eqnarray}

On the other hand, by \eqref{3.3} and the fact that
$p(r_0(\theta), \theta) = 0$, we have
$$p(r, \theta) > 0$$
for all $(r, \theta)$ with $r > r_0(\theta)$, which together with
\eqref{z5} result in that $p(\overline{r}, \overline{\theta}) >
0$. Since $(\overline{r}, \overline{\theta})$ is a point on the
bubble, we arrive at a contradiction.

Summing up, we in fact proved the following theorem.

\begin{thm}
The Riemann problem \eqref{1.1}, \eqref{a1} for the pressure
gradient equation admits a unique smooth solution. The pressure of
the solution is strictly positive in $\xi>0, \eta>0$.
\end{thm}

\section{Acknowledgement}
This work was done when Zhen Lei was visiting the Department of
Mathematics of Pennsylvania State University. He would like to
express his thanks for its hospitality. In particular, he thanks
Professor Chun Liu for his encouragement.
 Zhen Lei was partially supported by the National Science Foundation of
China under grant 10225102 and Foundation for Excellent Doctoral
Dissertation of China. Yuxi Zheng was partially supported by
National Science Foundation under grants NSF-DMS 0305479 and
NSF-DMS 0305114.


\end{document}